\documentclass{birkart}

\usepackage{amsmath}
\usepackage{amssymb}
\usepackage{amsfonts}

 \newtheorem{thm}{Theorem}[section]
 \newtheorem{cor}[thm]{Corollary}
 \newtheorem{lem}[thm]{Lemma}
 \newtheorem{prop}[thm]{Proposition}
 \theoremstyle{definition}

 \theoremstyle{remark}
 \newtheorem{rem}[thm]{Remark}

\numberwithin{equation}{section}
\numberwithin{figure}{section}

\newcommand{\D}{{\mathbb D}}

\newcommand{\T}{{\mathbb T}}
\newcommand{\R}{{\mathbb R}}

\newcommand{\re}{\operatorname{Re}}

\newcommand{\diff}{{\mathrm d}}
\newcommand{\imag}{{\mathrm i}}

\newcommand{\HSexp}{\operatorname{HSexp}}

\begin{document}
\title[Hele-Shaw flow]
{Hele-Shaw flow on weakly hyperbolic surfaces}
%----------Author 1
\author[Hedenmalm]
{H\aa{}kan Hedenmalm}

\address{Hedenmalm: Department of Mathematics\\
The Royal Institute of Technology\\
S -- 100 44 Stockholm\\
Sweden}

\email{haakanh@math.kth.se}

\thanks{Both authors wish to thank the G\"oran Gustafsson 
Foundation for generous support.}

%----------Author 2
\author[Olofsson]{Anders Olofsson}
\address{Olofsson: Department of Mathematics\\
The Royal Institute of Technology\\
S -- 100 44 Stockholm\\
Sweden}
\email{ao@math.kth.se}



%\date{March 23, 2004}
%\dedicatory{To my parents}
%%% ----------------------------------------------------------------------

\begin{abstract}
We consider the Hele-Shaw flow that arises from injection of two-dimensional
fluid into a point of a curved surface. The resulting fluid domains have
and are more or less determined implicitly by a mean value property for
harmonic functions. We improve on the results of Hedenmalm and Shimorin
\cite{HS} and obtain essentially the same conclusions while imposing a
weaker curvature condition on the surface. Incidentally, the curvature 
condition is the same as the one that appears in Hedenmalm and Perdomo's
paper \cite{HP}, where the problem of finding smooth area minimizing surfaces
for a given curvature form under a natural normalizing condition was 
considered. Probably there are deep reasons behind this coincidence.
\end{abstract}

%%% ----------------------------------------------------------------------
\maketitle
%%% ----------------------------------------------------------------------

\section{Introduction}

%\addtolength{\oddsidemargin}{-1.3cm}
%\addtolength{\evensidemargin}{-1.3cm}
%\addtolength{\textwidth}{1.5cm} 
%\addtolength{\textheight}{2.5cm}
%\addtolength{\topmargin}{-1.0cm}
%\addtolength{\footskip}{1.5cm}

%\begin{document}

%%% Research support

%\thanks{Research supported by the G\"oran Gustafsson Foundation.} 

%\subjclass{Primary: 31A30; Secondary: 30E10}

%\keywords{Hele-Shaw flow, weighted biharmonic Green function}

%\begin{abstract}
%\end{abstract}

%\maketitle

%%% Section numbering
%\setcounter{section}{-1}

%\section{Introduction}

%\setcounter{equation}{0}
%\setcounter{thm}{0}
%\setcounter{prop}{0}
%\setcounter{lemma}{0}
%\setcounter{cor}{0}
%\setcounter{remark}{0}

Let ${\mathbf\Omega}$ be a simply connected Riemann surface with a 
$C^\infty$-smooth metric $\diff{\mathbf s}$. 
We consider the case when the surface ${\mathbf\Omega}$ is conformally 
equivalent to a simply connected planar domain $\Omega$, which we pick to be
the unit disk $\Omega=\D$, and the metric $\diff{\mathbf s}$ is given by
\begin{equation}\label{metric}
\diff{\mathbf s}(z)^2=\omega(z)\,\vert \diff z\vert^2,\qquad z\in\Omega,
\end{equation}
for some positive $C^\infty$-smooth weight function $\omega$ in $\Omega$. 
This means that we use {\sl isothermal coordinates}, which is basically 
always possible;
the only real restriction imposed is that we do not allow ${\mathbf\Omega}$ to
be conformally  equivalent to the sphere or the plane. The {\sl area form} 
of the surface ${\mathbf\Omega}$ is (with suitable normalization)
$$\diff {\mathbf\Sigma}(z)=\omega(z)\,\diff\Sigma(z),\qquad z\in\Omega=\D,$$
where
$$\diff\Sigma(z)=\frac{\diff x\diff y}{\pi},\qquad z=x+\imag y,$$
is the normalized area element in the plane.  

The Gaussian curvature function associated to the metric $\diff{\bf s}$ 
given by (\ref{metric}) is the function ${\boldsymbol\kappa}$ defined by
$$
{\boldsymbol\kappa}(z)=
-\frac{2}{\omega(z)}\,\Delta(\log\omega)(z),\qquad z\in\D,
$$
where 
$$\Delta=\Delta_z=\frac{\partial^2}{\partial z\partial\bar z}
=\frac14\,\bigg(\frac{\partial^2}{\partial x^2}+
\frac{\partial^2}{\partial y^2}\bigg),\qquad z=x+\imag y,$$ 
is the (normalized) Laplacian in the plane. The curvature form associated 
to the metric $\diff{\bf s}$ given by (\ref{metric}) is the $2$-form 
${\mathbf K}$ defined by
$$
{\mathbf K}(z)={\boldsymbol\kappa}(z)\,\diff{\mathbf\Sigma}(z)=
-2\,\Delta(\log\omega)(z)\,\diff\Sigma(z),\qquad z\in\D.
$$

In the recent paper \cite{HP}, Hedenmalm and Perdomo studied the following 
problem: Suppose the curvature form $\mathbf K$ is given and 
$C^\infty$-smooth on the closed disk $\bar\D$, while the metric 
$\diff{\mathbf s}$ is to be found. If we ask of this metric which we are 
looking for that it is $C^\infty$-smooth on $\bar\D$ and has the 
{\sl mean value property}
\begin{equation}
h(0)=\int_\D h(z)\,\diff {\mathbf \Sigma}(z)=\int_\D h(z)\,\omega(z)\,
\diff {\Sigma}(z)
\label{eq-meanvalue}
\end{equation}
for all bounded harmonic functions $h$ in $\D$, can we then find such a 
metric, and is it unique? It turns out that this is indeed so if the condition
\begin{equation}\label{curvaturecond}
{\mathbf K}(z)+\frac{1}{2}\,{\mathbf K}_\mathbb{H}(z)\leq0,\qquad z\in\D,
\end{equation}
is met, where ${\mathbf K}_\mathbb{H}$ denotes the curvature form for 
the surface $\D$ when equipped with the (Poincar\'e) metric 
of constant curvature $-1$:
$${\mathbf K}_\mathbb{H}(z)=-\frac{4\,\diff\Sigma(z)}{(1-\vert z\vert^2)^2},
\qquad z\in\D.$$
Moreover, the result is fairly sharp, in the following sense: If the number 
$\frac12$ is replaced by the slightly larger number $0.52$ in condition 
(\ref{curvaturecond}), the resulting condition is too weak to force the 
existence of a $C^\infty$-smooth bordered surface $\mathbf\Omega=\langle
\D,\diff{\mathbf s}\rangle=\langle\D,\omega\rangle$ with the mean value 
property for harmonic functions and the given curvature form 
${\mathbf K}$. It was conjectured in \cite{HP} that the numerically 
obtained number $0.52$ could be replaced by $\frac12+\varepsilon$, for any 
positive $\varepsilon$.
We shall at times refer to (\ref{curvaturecond}) as the condition of 
{\sl weak hyperbolicity}. 

It is a natural question to ask whether there are subdomains $D(t)$, depending
on the parameter $t$, $0<t<1$, with the modified mean value property
\begin{equation}
t\,h(0)=\int_{D(t)} h(z)\,d{\mathbf \Sigma}(z)=\int_{D(t)} h(z)\,\omega(z)\,
d{\Sigma}(z)
\label{eq-meanvalue-t}
\end{equation}
for all bounded harmonic functions $h$ in $D(t)$. We make the {\sl a priori}
assumption that the domains $D(t)$ are simply connected with $C^\infty$-smooth
boundaries. Then, if we look at the curvature condition (\ref{curvaturecond}) 
and use the fact that the curvature form for the Poincar\'e metric decreases 
as the domain gets smaller, we see that the analogous curvature condition 
for ${\mathbf\Omega}(t)=\langle D(t),\omega\rangle$ is fulfilled, so by the 
Hedenmalm-Perdomo theorem, the surfaces ${\mathbf\Omega}(t)$ exist 
as abstract surfaces, being determined uniquely by the given curvature form 
and the mean value property (\ref{eq-meanvalue-t}). 
However, this does not tell us that the abstract surface 
${\mathbf\Omega}(t)$ forms a subregion of $\mathbf\Omega$. It would be 
desirable to know that this is so. 
If this works, then perhaps it will eventually be possible to understand the 
Hedenmalm-Perdomo theorem in terms of the following growth process: 

\noindent (1) for $t$ close to $0$, we form the surface ${\mathbf\Omega}(t)$, 
by fairly elementary means; 

\noindent (2) as $t$ increases, we iteratively add abstractly infinitesimally
thin boundary layers to the surface ${\mathbf\Omega}(t)$, while checking 
that the curvature condition (\ref{curvaturecond}) permits us to do so.

Here, we shall obtain the existence of the smooth simply connected subdomains 
$D(t)$ of $\D$ with (\ref{eq-meanvalue-t}) under the weak hyperbolicity 
condition (\ref{curvaturecond}) 
in the category of {\sl real-analytic surfaces}. Moreover, we shall prove 
that the boundary curves $\partial D(t)$ are real-analytically smooth, and
that we get an extra degree of analytic smoothness at the origin for 
$t\to0$. The growth of the domains $D(t)$ is a certain kind of {\sl Hele-Shaw
flow}, modelling the growth of a patch of two-dimensional fluid as fluid
is injected at the origin. We extend the Hele-Shaw exponential mapping theorem
obtained by Hedenmalm and Shimorin in \cite{HS} to weakly hyperbolic surfaces.
From the point of view of differential geometry, it is perhaps surprising that
this is possible, because it does not work globally for the ordinary 
exponential mapping. After all, it is rather easy to construct a 
real-analytically smooth surface ${\mathbf\Omega}=\langle\D,\omega\rangle$ 
whose curvature form $\mathbf K$ satisfies (\ref{curvaturecond}) while two 
geodesics emanating from the origin meet at some other point; we supply an 
example of such a situation in Section~\ref{exwhs}. Probably this discrepancy
can be explained by the viscosity of the fluid as compared with the lack
thereof in the photonic gas, which is the physical manifestation of the metric
flow. 

The techniques used in the paper are basically similar to what was used
by Hedenmalm-Shimorin in \cite{HS}. The main difference is that the Green 
function $\Gamma_1$ for the weighted biharmonic operator $\Delta(1-|z|^2)^{-1}
\Delta$ is used in place of the biharmonic Green function $\Gamma$. We
suspect that the main results of this paper are false if the constant 
$\frac12$ in the condition (\ref{curvaturecond}) is replaced by 
$\frac12+\varepsilon$, for any fixed positive $\varepsilon$. In another vein,
we suspect that the results obtained here for real-analytic surfaces will 
remain valid for $C^\infty$-smooth surfaces as well.  

\section{Preliminaries on Hele-Shaw domains}\label{prelHSdomains}

Let $\omega$ be a $C^\infty$-smooth weight function on $\bar\D$ which is 
strictly positive (that is, it is positive at all points).
We denote by $G$ the Green function for the Laplacian $\Delta$ on $\D$:
\begin{equation*}
G(z,\zeta)=\log\left\vert\frac{z-\zeta}{1-\bar\zeta z}\right\vert^2,
\qquad(z,\zeta)\in\bar\D\times\D.
\label{Greenfn}
\end{equation*}
For $0<t<+\infty$, we consider the function
$$
V_t(z)=t\,G(z,0)-
\int_\D G(z,\zeta)\,\omega(\zeta)\,\diff\Sigma(\zeta),
\qquad z\in\D.
$$
This function $V_t$ is smooth in $\D\setminus\{0\}$, vanishes on 
the unit circle $\T=\partial\D$, and has a logarithmic singularity at the 
origin.
We denote by $\widehat{V}_t$ the smallest superharmonic majorant of $V_t$ 
in $\D$. The set $D(t)$ is now defined by 
\begin{equation*}%\label{dfnhsdomain}
D(t)=D(t;\omega)=
\Big\{z\in\D:\,\, V_t(z)<\widehat{V}_t(z)\Big\}.
\end{equation*}
The sets $D(t)$ are known to have a number of basic 
properties (see~\cite{HS}). 
For instance, the sets $D(t)$ are open, connected and increasing in the 
positive parameter  $t$. 

Associated to the domains $D(t)$, we consider the number
$$
T=T(\omega)=\sup\Big\{t\in]0,+\infty[:\  D(t)\Subset\D\Big\},
$$
the {\sl termination time} for the Hele-Shaw flow. It is known that 
$0<T\leq+\infty$ (see~\cite[Proposition~2.8]{HS}).
For $0<t<T$, the domain $D(t)$ is called a {\sl Hele-Shaw domain}. 
For $t\geq T$, the domain $D(t)$ is sometimes referred to as a 
{\sl generalized Hele-Shaw domain}. 
%We shall later consider the case when $\D$ has the mean-value 
%property (\ref{eq-meanvalue}), which implies that $0<T\le1$. 

\section{A weighted biharmonic Green function}\label{wbhGfcn}

In this section we shall need some properties of the Green function 
for the weighted biharmonic operator 
$\Delta(1-\vert z\vert^2)^{-1}\Delta$. 
We review the results needed and refer to~\cite{HP} for details.

The Green function for the weighted biharmonic operator 
$\Delta(1-\vert z\vert^2)^{-1}\Delta$ is the 
function $\Gamma_1$ on $\bar\D\times\D$ 
solving (in an appropriate sense) for fixed $\zeta\in\D$ 
the boundary value problem
\begin{equation}
\left\{
\begin{array}{rclr}
\Delta_z(1-\vert z\vert^2)^{-1}\Delta_z\Gamma_1(z,\zeta)&=&
\delta_\zeta(z)&\quad z\in\D,\\
\Gamma_1(z,\zeta)&=&0 & z\in\T,\\
\nabla_z\Gamma_1(z,\zeta)&=&0 & z\in\T;
\end{array}
\right.
\end{equation}
here $\delta_\zeta$ is the (unit) Dirac mass at $\zeta$ and 
$\T=\partial\D$ is the unit circle.
The Green function $\Gamma_1$ has the explicit expression
\begin{align}\notag
\Gamma_1(z,\zeta)&=
\bigg\{\vert z-\zeta\vert^2-
\frac{1}{4}\,\big\vert z^2-\zeta^2\big\vert^2\bigg\}G(z,\zeta)+
\frac{1}{8}\big(1-\vert z\vert^2\big)\big(1-\vert\zeta\vert^2\big)\\
&\times\bigg\{7-\vert z\vert^2-\vert\zeta\vert^2-\vert z\zeta\vert^2-
4\re(z\bar\zeta)-2\big(1-\vert z\vert^2\big)\big(1-\vert\zeta\vert^2\big)
\frac{1-\vert z\zeta\vert^2}{\vert1-z\bar\zeta\vert^2}\bigg\},
\label{G1dfn}
\end{align}
where 
$$
G(z,\zeta)=\log\left\vert\frac{z-\zeta}{1-\bar\zeta z}\right\vert^2,
\qquad(z,\zeta)\in\bar\D\times\D,
$$ 
is the Green function for the Laplacian in $\D$. 
We mention that the above 
explicit expression for $\Gamma_1$ first appeared in~\cite{Hprobl}.
A derivation of formula (\ref{G1dfn}) can be found in~\cite{W2}.

It is known that the Green function $\Gamma_1$ is positive in the bidisk:
$$
\Gamma_1(z,\zeta)>0,\qquad(z,\zeta)\in\D\times\D;
$$
see \cite[Proposition~3.4]{Hprobl} or \cite[Lemma~2.2]{HP}.
%We shall also need the growth estimate 
%\begin{equation}\label{G1growth}
%\Gamma_1(z,\zeta)=O\big((1-\vert z\vert)^3\big),
%\quad\vert z\vert\to1,
%\end{equation}
%where $\zeta\in\D$ is fixed.

The function $H_1$ defined by
$$
H_1(z,\zeta)=\big(1-\vert\zeta\vert^2\big)
\bigg\{\frac{1}{2}\big(3-\vert\zeta\vert^2\big)
\frac{1-\vert z\zeta\vert^2}{\vert1-z\bar\zeta\vert^2}+
\big(1-\vert\zeta\vert^2\big)\re\bigg[\frac{z\bar\zeta}
{(1-z\bar\zeta)^2}\bigg]\bigg\},
$$
for $(z,\zeta)\in\bar\D\times\bar\D$ with $z\neq\zeta$, has been coined the 
{\sl harmonic compensator}; it is positive on $\bar\D\times\bar\D$. 
Indeed, by harmonicity in $z$, it is enough to consider 
$(z,\zeta)\in\T\times\D$, in which case we have that
\begin{align*}
H_1(z,\zeta)&\geq\big(1-\vert\zeta\vert^2\big)
\bigg\{\frac{1}{2}\big(3-\vert\zeta\vert^2\big)
\frac{1-\vert\zeta\vert^2}{\vert 1-z\bar\zeta\vert^2}-
\big(1-\vert\zeta\vert^2\big)
\frac{\vert\zeta\vert}{\vert 1-z\bar\zeta\vert^2}\bigg\}=\\
&=\frac{1}{2}\,\frac{(1-\vert\zeta\vert^2)^2}{\vert1-z\bar\zeta\vert^2}
\big(1-\vert\zeta\vert\big)\big(3+\vert\zeta\vert\big)>0.
\end{align*}
The functions $\Gamma_1$ and $H_1$ are related by the formula
\begin{equation}\label{lplG1}
\Delta_z\Gamma_1(z,\zeta)=
\big(1-\vert z\vert^2\big)\big(G(z,\zeta)+H_1(z,\zeta)\big),
\qquad(z,\zeta)\in\bar\D\times\D,
\end{equation}
which can be verified by straightforward computation (see~\cite{HP}).

The following lemma establishes an integral representation 
in terms of the functions $\Gamma_1$ and $H_1$.

\begin{lem}\label{reprformula}
Let $u$ be a smooth function on $\bar\D$ such that $u=0$ on $\T$. 
Then $u$ admits the representation
$$
u(\zeta)=\int_\D\Gamma_1(z,\zeta)\,\Delta(1-\vert z\vert^2)^{-1}
\Delta u(z)\,\diff\Sigma(z)+\frac{1}{2}\int_\T H_1(z,\zeta)\,
\partial_n u(z)\,\diff\sigma(z),
\quad\zeta\in\D,
$$
where $\partial_n$ denotes differentiation in the 
inward normal direction and $\diff\sigma$ is normalized arc length measure on 
$\T$.
\end{lem}

A more general representation formula is planned to appear elsewhere.
Therefore we omit the proof of Lemma~\ref{reprformula}.
The lemma can also be obtained by repeated applications 
of Green's formula.
We recall that with our normalizations Green's formula
(Green's second identity) takes the form
\begin{equation}\label{Greenformula}
\int_\Omega\big(u\Delta v-v\Delta u\big)\,\diff\Sigma=
\frac{1}{2}\int_{\partial\Omega}\big(v\partial_n u-
u\partial_nv\big)\,\diff\sigma,
\end{equation}
where $\diff\sigma=\vert\diff z\vert/2\pi$.

The following lemma is a simple consequence of the fact that $\Gamma_1$ is
positive. The technique is analogous to that which was developed for the
biharmonic Green function $\Gamma$ in \cite{DKSS}; see also \cite{H91,H92,HP}.

\begin{lem}\label{nucontremb}
Let $\nu$ be a smooth function on $\bar\D$ such that 
$$z\mapsto\frac{\nu(z)}{1-\vert z\vert^2}$$
is subharmonic in $\D$.
Assume also that $\nu$ is reproducing at the origin in the sense that
$$
h(0)=\int_\D h(z)\,\nu(z)\,\diff\Sigma(z)
$$
holds for every harmonic polynomial $h$.
Then the inequality 
$$
\int_\D u(z)\,2\big(1-\vert z\vert^2\big)\,\diff\Sigma(z)\leq
\int_\D u(z)\,\nu(z)\,\diff\Sigma(z)
$$
holds for every $u\in C^\infty(\bar\D)$ which is subharmonic in $\D$.
\end{lem}

\begin{proof}
We consider the function $\Phi$ defined by 
$$
\Phi(z)=
\int_\D G(z,\zeta)\,\big[\nu(\zeta)-2\big(1-\vert\zeta\vert^2\big)\big]\,
\diff\Sigma(\zeta),
\qquad z\in\D,
$$
where $G$ is the Green function for the Laplacian in $\D$.
By construction $\Phi\vert_\T=0$. 
Since $\Delta\Phi=\nu-2(1-\vert z\vert^2)$ annihilates harmonic polynomials, 
an application of Green's formula (\ref{Greenformula}) 
shows that also $\partial_n\Phi\vert_\T=0$.
By Lemma~\ref{reprformula}, the function $\Phi$ has the representation
\begin{multline*}
\Phi(z)=\int_\D\Gamma_1(z,\zeta)\,
\Delta_\zeta(1-\vert\zeta\vert^2)^{-1}
\Delta_\zeta\Phi(\zeta)\,\diff\Sigma(\zeta)\\
=\int_\D\Gamma_1(z,\zeta)\,
\Delta_\zeta\bigg[\frac{\nu(\zeta)}{1-\vert\zeta\vert^2}\bigg]\,
\diff\Sigma(\zeta)\geq0.
\end{multline*}
Another application of Green's formula now shows that
\begin{multline*}
\int_\D u(z)\big[\nu(z)-2\big(1-\vert z\vert^2\big)\big]\,\diff\Sigma(z)=
\int_\D u(z)\,\Delta\Phi(z)\,\diff\Sigma(z)\\=
\int_\D\Delta u(z)\,\Phi(z)\,\diff\Sigma(z)\geq0,
\end{multline*}
which concludes the proof.
\end{proof}

In the proof of the next lemma, we shall use some easy properties 
of Bergman kernel functions. Let us denote by $A^2_\alpha(\D)$, for 
$-1<\alpha<\infty$, the weighted Bergman space of all analytic functions in 
$\D$ that are square integrable with respect to the measure
$$\diff \Sigma_\alpha(z)=\omega_\alpha(z)\,\diff \Sigma(z),$$
where the weight $\omega_\alpha$ is as follows:
$$
\omega_\alpha(z)=(1+\alpha)\big(1-\vert z\vert^2\big)^\alpha,\qquad z\in\D.
$$ 
The Bergman kernel function for the space $A^2_\alpha(\D)$ is the function
$$
K_\alpha(z,\zeta)=\frac{1}{\big(1-z\bar\zeta\big)^{2+\alpha}},
\qquad(z,\zeta)\in\D\times\D;
$$
see~\cite[Proposition~1.1.4]{HKZ}. The reproducing property of the kernel 
function $K_\alpha$ asserts that 
$$
f(z)=\int_\D K_\alpha(z,\zeta)\,f(\zeta)\,\omega_\alpha(z)\,
\diff\Sigma(\zeta),\qquad z\in\D,
$$
for $f\in A^2_\alpha(\D)$.
%We can now prove the following lemma.

The following lemma is inspired by~\cite[Proposition~3.2]{HP}.

\begin{lem}\label{nuonT}
Let $\nu$ be a nonnegative smooth weight function satisfying 
the assumptions of Lemma~\ref{nucontremb}.
Let $e^{\imag\theta}\in\T$.
Then either $\nu(e^{\imag\theta})>0$, or 
$\nu(e^{\imag\theta})=0$ and $\partial_n\nu(e^{\imag\theta})>0$.
\end{lem}

\begin{proof}
Without loss of generality, we may assume that $e^{\imag\theta}=1$. To
reach a contradiction, we assume that the conclusion of the lemma does not 
hold. Then there is a positive constant $C$ such that 
$$\nu(z)\leq C\vert z-1\vert^2,\qquad z\in\D.$$

Fix $r$ in the interval $0<r<1$. By Lemma~\ref{nucontremb} applied to the 
function 
$$u(z)=\big\vert K_1(z,r)\big\vert^2=\frac{1}{|1-rz|^6},$$ 
we have that
\begin{equation}\label{absurdineq}
\frac1{(1-r^2)^3}=\int_\D\big\vert K_1(z,r)\big\vert^2\,2(1-\vert z\vert^2)
\,\diff\Sigma(z)\leq 
C\int_\D\big\vert K_1(z,r)\big\vert^2\,\vert z-1\vert^2\, \diff\Sigma(z).
\end{equation}
The leftmost equality in (\ref{absurdineq}) holds because of the 
reproducing property of $K_1$. We turn to estimating the right hand side 
in (\ref{absurdineq}). By the elementary inequality 
$$r^2\vert z-1\vert^2\leq\vert 1-rz\vert^2,\qquad z\in\D,$$ 
we have that
$$
\int_\D\big\vert K_1(z,r)\big\vert^2\,\vert z-1\vert^2\, \diff\Sigma(z)\leq
\frac{1}{r^2}\int_\D\frac{1}{\vert 1-rz\vert^4}\,
\diff\Sigma(z)=\frac{1}{r^2(1-r^2)^2},
$$
where the last equality follows from the reproducing property of $K_0$.
Thus, by (\ref{absurdineq}), we have that 
$$
\frac{1}{(1-r^2)^3}\leq\frac{C}{r^2(1-r^2)^2},
$$
which is clearly impossible for $r$ close to $1$.
\end{proof}

The next lemma is analogous to \cite[Lemma~5.2]{HS}.

\begin{lem}\label{wbhgestimate}
Let $u$ be a smooth real-valued function in $\bar\D\setminus\{0\}$ 
with a logarithmic singularity at the origin such that  
$$
\Delta(1-\vert z\vert^2)^{-1}\Delta u=\Delta\delta_0-\mu\qquad\text{in}\ \D,
$$ 
where $0\leq\mu\in C^\infty(\D)$. 
Suppose that $u\vert_\T=0$ and $\partial_n u\leq0$ on $\T$.
Then
$$
u(z)\leq\log\vert z\vert^2+\frac{3}{2}-2\vert z\vert^2+
\frac{1}{2}\vert z\vert^4<0,\qquad z\in\D.
$$
\end{lem}

\begin{proof}
Recall that the kernels $\Gamma_1$ and $H_1$  are both positive.
By the representation formula in Lemma~\ref{reprformula} we have that
$$
u(\zeta)\leq\int_\D\Gamma_1(z,\zeta)\,\Delta\delta_0(z)\,\diff\Sigma(z)=
\Delta_z\Gamma_1(0,\zeta)=\log\vert\zeta\vert^2+\frac{3}{2}-
2\vert\zeta\vert^2+\frac{1}{2}\vert\zeta\vert^4,
$$ 
where the integration is to be interpreted in the usual sense of 
distribution theory. The proof is complete.
\end{proof}

\section{Simply connectedness of Hele-Shaw domains}\label{simplyconnHSdomain}

In this section we shall show that each Hele-Shaw domain $D(t)$, $0<t<T$, 
is simply connected and has real-analytic Jordan boundary.
First we need a lemma.

\begin{lem}\label{Westimate}
Let $\nu$ be a strictly positive smooth function on $\bar\D$ such that 
the function 
$$z\mapsto\frac{\nu(z)}{1-\vert z\vert^2}$$ 
is subharmonic in $\D$. Let 
$$
W(z)=\log\vert z\vert^2-\int_\D G(z,\zeta)\,\nu(\zeta)\,
\diff\Sigma(\zeta),
\qquad z\in\D,
$$ 
and denote by $\widehat{W}$ the smallest superharmonic majorant 
of $W$ in $\D$.
Assume that the coincidence set
$$\Big\{z\in\D:\ W(z)=\widehat{W}(z)\Big\}$$ 
is a compact subset of $\D$. Then 
$$W(z)\le\log\vert z\vert^2+\frac{3}{2}-2\vert z\vert^2+
\frac{1}{2}\vert z\vert^4<0,\qquad z\in\D;$$
in particular, the coincidence set is empty.
\end{lem}

\begin{proof}
We first introduce some preliminary notation. 
We consider the functions 
$$
W_r(z)=r\log\vert z\vert^2-\int_\D G(z,\zeta)\,\nu(\zeta)\,\diff\Sigma(\zeta),
\qquad z\in\D,
$$
for $1\leq r<\infty$, and note that $W_1=W$.
By construction, $W_r$ vanishes on $\T$.
Denote by $\widehat{W}_r$ the smallest superharmonic majorant of 
$W_r$ in $\D$, and write 
$$B(r)=\Big\{z\in\D:\ W_r(z)<\widehat{W}_r(z)\Big\}.$$
The set $B(r)$ gets bigger as $r$ increases; see~\cite[Proposition~2.7]{HS}.
If $W_r\leq0$ in $\D$, then $\partial_n W_r\leq0$ on $\T$.
Conversely, if  $\partial_n W_r\leq0$ on $\T$, 
then Lemma~\ref{wbhgestimate} applied to the function $r^{-1}W_r$ 
shows that 
\begin{equation}\label{Wrestimate}
W_r(z)\leq r\Big(\log\vert z\vert^2+\frac{3}{2}-2\vert z\vert^2+
\frac{1}{2}\vert z\vert^4\Big)<0,\qquad z\in\D.
\end{equation}

If $\partial_n W_1\leq0$ on $\T$, we are done. 
Assume, to reach a contradiction, that 
$\max_\T\partial_n W_1$ is strictly positive.
In view of the definition of $W_r$, we have that 
$$
\partial_n W_r(z)=-2(r-1)+\partial_n W_1(z),\qquad z\in\T.
$$
This formula makes evident that 
$$\max_\T\partial_n W_r>0\quad \text{for}\quad
1\leq r<r_1=1+\max_\T\partial_n W_1/2$$ 
while $\max_\T\partial_n W_r\leq0$ for $r\geq r_1$.

For $1\leq r<r_1$, the function $W_r$ attains a positive maximum at some 
point $z(r)\in\D$. Clearly, $z(r)\in\D\setminus B(r)\subset \D\setminus B(1)$.
Let $\{r_j\}_{j=2}^\infty$ be a sequence with $1\leq r_j<r_1$ and
$r_j\to r_1$ as $j\to+\infty$. In view of the assumed compactness of 
$\D\setminus B(1)$, we may by passing to a subsequence assume that 
$z(r_j)$ converges to a point $z_1\in\D\setminus B(1)$ as $j\to+\infty$. 
Since $W_{r_j}(z(r_j))>0$, we obtain in the limit that $W_{r_1}(z_1)\geq0$, 
which contradicts (\ref{Wrestimate}) for $r=r_1$. It follows that $B(1)=\D$,
and that the estimate (\ref{Wrestimate}) holds for all $r$, $1\le r<+\infty$.
\end{proof}

We may now derive the asserted properties of $D(t)$. We recall the definition 
of the termination time $T=T(\omega)$ for the Hele-Shaw flow.

\begin{thm}\label{simplyconnected}
Let ${\mathbf\Omega}=\langle\D,\omega\rangle$ be a 
simply connected Riemann surface with a metric 
$\diff{\mathbf s}^2=\omega\vert\diff z\vert^2$. 
%given by (\ref{metric}).
Assume that the weight function $\omega$ is real-analytic and 
strictly positive in $\D$ and that the curvature condition 
\begin{equation}\label{curvaturecondition}
{\mathbf K}(z)+\frac{1}{2}\,{\mathbf K}_\mathbb{H}(z)\leq0,
\qquad z\in\D,
\end{equation}
is satisfied.
Then, for $0<t<T$, 
the Hele-Shaw domain $D(t)=D(t;\omega)$ is simply connected.
\end{thm}

\begin{proof}
Without loss of generality, we can assume that $\omega$ is smooth up to the 
boundary. It is known that $\partial D(t)$ has a local Schwarz function 
at every point and that this implies that $D(t)$ 
can have at most finitely many holes 
(see~\cite[Section~4]{HS}; the argument uses the work of Sakai~\cite{Sakai}).
Let $D_{\bullet}(t)$ be the simply connected domain obtained from 
$D(t)$ by adding all the interior holes. 
Let $\varphi:\D\to D_\bullet(t)$ be a Riemann map with $\varphi(0)=0$. 
In view of the regularity of $\partial D(t)$ implied by the existence of a
local Schwarz function, the map $\varphi$ extends analytically to a 
neighborhood of the closed disk $\bar\D$. Let $B=\varphi^{-1}(D(t))$, 
and note that $\D\setminus B$ is a compact subset of $\D$. 
Introduce the function  
\begin{equation}\label{dfnnu}
\nu(z)=\frac{1}{t}\,
(\omega\circ\varphi)(z)\,\vert\varphi'(z)\vert^2,\qquad z\in\D.
\end{equation}
By the curvature assumption (\ref{curvaturecondition}), 
the function 
$$z\mapsto\frac{\nu(z)}{1-\vert z\vert^2}$$ 
is logarithmically subharmonic in $\D$ and, hence, subharmonic there. 

The set $B$ can be interpreted as the non-coincidence set of 
an obstacle problem. A computation shows that 
$$
\Delta\big[V_t\circ\varphi(z)\big]=t\,\delta_0(z)-
\omega\circ\varphi(z)\,\,\vert\varphi'(z)\vert^2,\qquad z\in \D,
$$
where the potential function $V_t$ is as in Section~\ref{prelHSdomains}.
Let $W$ be defined in terms of $\nu$ as in Lemma~\ref{Westimate}.
It is easy to see that
$$
W(z)=\frac{1}{t}
\Big(V_t\circ\varphi(z)-P\big[V_t\circ\varphi\vert_\T\big](z)\Big),
\qquad z\in\D.
$$ 
where $P[\cdot]$ denotes the usual Poisson integral in $\D$. 
By~\cite[Proposition~2.9(a)]{HS}, the function $\widehat{V}_t$ 
is also the smallest superharmonic majorant for $V_t$ in $D_\bullet(t)$.
By conformal invariance, the function $\widehat{V}_t\circ\varphi$ 
is the smallest superharmonic majorant for $V_t\circ\varphi$ in $\D$, 
and we have that 
$$
\widehat W(z)=\frac{1}{t}
\Big(\widehat V_t\circ\varphi(z)-P\big[V_t\circ\varphi\vert_\T\big](z)\Big),
\qquad z\in\D.
$$ 
It is now clear that 
$$B=\Big\{z\in\D:\ W(z)<\widehat W(z)\Big\}.$$ 
Lemma~\ref{Westimate} shows that $B=\D$, which means that 
$D(t)=D_\bullet(t)$ is simply connected. The proof is complete.
\end{proof}

\begin{thm}
Let ${\mathbf\Omega}=\langle\D,\omega\rangle$ be as in  
Theorem~\ref{simplyconnected}.
Then the boundary $\partial D(t)$ of the Hele-Shaw domain is a real-analytic 
Jordan curve for $0<t<T$.
\label{realanalytic}
\end{thm}

\begin{proof}
In view of the apriori regularity of $\partial D(t)$ 
(see~\cite[Section~4]{HS}) which follows from Sakai's theorem \cite{Sakai}, 
as well as from Theorem~\ref{simplyconnected}, the boundary $\partial D(t)$ 
is real-analytic with the exception of at most finitely many cusp 
or contact points (see \cite{HS} for the terminology). 
Let $\varphi:\D\to D(t)$ be a Riemann map with $\varphi(0)=0$. 
The regularity of $\partial D(t)$ afforded by the existence of a local
Schwarz function shows that the map $\varphi$ extends analytically to a 
neighborhood of $\bar\D$. The cusp points of $\partial D(t)$ correspond to 
points $z\in\T$ such that $\varphi'(z)=0$.
Arguing as in~\cite[Subsection~5.3]{HS} 
we see that the mean value identity 
$$
h(0)=\int_\D h(z)\,\nu(z)\,\diff\Sigma(z)
$$
holds for harmonic polynomials $h$, where $\nu$ is given by (\ref{dfnnu}).
The non-vanishing of $\varphi'$ on $\T$ now follows by Lemma~\ref{nuonT}.

If $D(t)$ has a contact point, then, for $t'=t+\delta$, 
with $\delta>0$ sufficiently small, the domain $D(t')$ has a hole, 
which is impossible by Theorem~\ref{simplyconnected}. For details of this
argument, see~\cite[Subsection~5.3]{HS}. 
\end{proof}

\section{Modification of a lemma of Korenblum}\label{Korenblum}

The following lemma is inspired by~\cite[Lemma~8.3]{HS} 
which is due to Boris Korenblum.

\begin{lem}\label{Klemma}
Fix $\alpha$, $0<\alpha<+\infty$. Let $\omega$ be a nonnegative function 
such that the function 
\begin{equation}\label{omegaash}
z\mapsto\log\frac{\omega(z)}{(1-\vert z\vert^2)^{2\alpha}}
\end{equation}
is subharmonic in $\D$. 
Assume also that $\omega$ is area-integrable in $\D$ and 
reproducing at the origin in the sense that
$$
h(0)=\int_\D h(z)\,\omega(z)\,\diff\Sigma(z)
$$
holds for all harmonic polynomials $h$.
Then, for a specific positive constant $c_{p,\alpha}$ depending only on 
$p,\alpha$, we have 
\begin{equation}\label{omegalengthestimate}
\int_0^1\omega(r)^p\, \diff r\leq c_{p,\alpha}<+\infty
\end{equation}
for $$0<p<\frac{\pi}{\pi+2\alpha(4-\pi)}.$$
\end{lem} 

\begin{rem} The constant $c_{p,\alpha}$ is given by equation (\ref{cpa})
below.
\end{rem}

\begin{proof}
The proof depends on the formula 
\begin{equation}\label{Pkformula}
\int_0^1\frac{1}{(1-r)^2}\,\chi_{\D(r,1-r)}(z)\,\diff r=
\frac{1-\vert z\vert^2}{\vert 1-z\vert^2},\qquad z\in\D,
\end{equation}
where $\chi_D$ denotes the characteristic function for the set $D$ and 
$\D(z_0,r)$ is the open disk with center $z_0$ and radius $r$. 
The formula (\ref{Pkformula}) is obtained by straightforward computation.

By the reproducing property of $\omega$, we have 
$$\int_\D\frac{1-\vert rz\vert^2}{\vert 1-rz\vert^2}
\,\omega(z)\,\diff\Sigma(z)=1,
\qquad 0<r<1,$$
and an application of Fatou's lemma shows that
\begin{equation}\label{omegabalayage}
\int_\D 
\frac{1-\vert z\vert^2}{\vert 1-z\vert^2}\,
\omega(z)\,\diff\Sigma(z)\leq1.
\end{equation}

We now turn to the proof of (\ref{omegalengthestimate}).
Fix an $r$ with $0<r<1$. By the subharmonicity assumption (\ref{omegaash}),
we have that
\begin{multline*}
\log\frac{\omega(r)}{(1-r^2)^{2\alpha}}\leq
\frac{1}{(1-r)^2}\int_{\D(r,1-r)}\log\omega(z)\,\diff\Sigma(z)\\
+\frac{2\alpha}{(1-r)^2}\int_{\D(r,1-r)}\log\frac{1}{1-\vert z\vert^2}\,
\diff\Sigma(z).
\end{multline*}
By an application of the geometric-arithmetic mean inequality,
we see that
\begin{equation}\label{omegaestimate}
\frac{\omega(r)}{(1-r^2)^{2\alpha}}\leq
%\exp\Big(\frac{1}{(1-r)^2}\int_{\D(r,1-r)}\log\frac{1}{1-\vert z\vert^2}
%\diff\Sigma(z)\Big)
F(r)^{2\alpha}
\frac{1}{(1-r)^2}\int_{\D(r,1-r)}\omega(z)\,\diff\Sigma(z),
\end{equation}
where we have written $F(r)$ for the quantity
\begin{equation}\label{Fdfn}
F(r)=\exp\bigg\{
\frac{1}{(1-r)^2}\int_{\D(r,1-r)}\log\frac{1}{1-\vert z\vert^2}\,
\diff\Sigma(z)\bigg\},
\qquad 0<r<1.
\end{equation}
It is clear that $F(r)$ is bounded from above away from the right end-point 
$r=1$. Below, we will show that $F(r)$ admits the estimate 
\begin{equation}\label{Festimate}
F(r)\leq \frac{1}{(1-r)^{4/\pi+\varepsilon}},\qquad r_\varepsilon<r<1,
\end{equation}
for every positive $\varepsilon$, where $0<r_\varepsilon<1$.
In view of (\ref{Pkformula}), (\ref{omegabalayage}), and 
(\ref{omegaestimate}), we have that
$$
\int_0^1\frac{\omega(r)}{(1-r^2)^{2\alpha}F(r)^{2\alpha}}\,\diff r\leq
\int_0^1\frac{1}{(1-r)^2}\int_{\D(r,1-r)}\omega(z)\,\diff\Sigma(z)
\,\diff r\leq1.
$$
An application of H\"older's inequality now gives that 
\begin{equation}
\int_0^1\omega(r)^p\, \diff r\leq
\bigg\{
\int_0^1\big[(1-r^2)F(r)\big]^{2\alpha p/(1-p)}\,\diff r
\bigg\}^{1-p}
=c_{p,\alpha},
\label{cpa}
\end{equation}
where the rightmost equality is used to define the constant $c_{p,\alpha}$.
The estimate (\ref{omegalengthestimate}) now follows. 
%with $c_{p,\alpha}=\int_0^1\big((1-r^2)F_\alpha(r)\big)^{p/(1-p)}\diff r$.
The finiteness of $c_{p,\alpha}$ %this last integral 
for the asserted values of $p$ 
%(for $0<p<\min(1,\pi/(8\alpha))$) 
follows by (\ref{Festimate}).

Finally, we turn to the proof of (\ref{Festimate}).
By a change to polar coordinates in (\ref{Fdfn}), we arrive at the formula
$$
\log F(r)=
\frac{2}{\pi(1-r)^2}\int_{2r-1}^1 t\log\bigg[\frac{1}{1-t^2}\bigg]
\arcsin\bigg[\frac{\sqrt{(1-t^2)(t^2-(2r-1)^2)}}{2rt}\bigg]\,\diff t
$$
which is valid for $\frac12<r<1$.
We note that by the geometric-arithmetic mean inequality, we have that
$$
\frac{\sqrt{(1-t^2)[t^2-(2r-1)^2]}}{2rt}\leq\frac{1-r}{t}.
$$
From an elementary estimate of the arcsine, we find that
$$
\log F(r)\leq\frac{2}{\pi}\,\frac{1+\varepsilon}{1-r}
\int_{2r-1}^1\log\bigg[\frac{1}{1-t^2}\bigg]\,\diff t
$$
for $r$, $\frac12<r<1$, sufficiently close to $1$. 
A straightforward computation shows that
$$
\int_{2r-1}^1\log\bigg[\frac{1}{1-t^2}\bigg]\,\diff t=
2(1-r)\log\frac{1}{1-r}+4(1-\log2)(1-r)+2r\log r.
$$
Using the elementary inequality $\log r\leq r-1$, we see that 
$$
\int_{2r-1}^1\log\bigg[\frac{1}{1-t^2}\bigg]\,\diff t\leq
2(1-r)\log\frac{1}{1-r}
$$
for $2(1-\log2)\leq r<1$.
As we return to the function $F$, we see that
$$
\log F(r)\leq\frac{4}{\pi}(1+\varepsilon)\log\frac{1}{1-r}
$$
for $r$ close to $1$, which yields (\ref{Festimate}). 
\end{proof}

\begin{rem} 
Assume that $\omega$ is logarithmically subharmonic in $\D$ and 
reproducing is the sense of Lemma~\ref{Klemma}. 
Letting $\alpha\to0$ in (\ref{omegalengthestimate}), we see that 
$$\int_0^1\omega(r)^p\diff r\leq1,\qquad 0<p<1.$$
By monotone convergence, we conclude that
$$\int_0^1\omega(r)\diff r\leq1.$$
In~\cite[Lemma~8.3]{HS}, this last inequality was shown to hold true 
under the slightly weaker assumption that 
$\omega$ is subharmonic and reproducing in $\D$. 
\end{rem}

In terms of the curvature form ${\mathbf K}$ of 
${\mathbf\Omega}=\langle\D,\omega\rangle$, the assumption (\ref{omegaash}) 
means that 
\begin{equation}\label{curvaturealpha}
{\mathbf K}(z)+\alpha\,{\mathbf K}_\mathbb{H}(z)\leq0,\qquad z\in\D,
\end{equation}
as is easily verified by straightforward computation.
In this context, we recall that the Poincar\'e metric of constant curvature 
$-1$ on $\D$ corresponds to the weight function 
$$
\omega_\mathbb{H}(z)=\frac{4}{(1-\vert z\vert^2)^2},\qquad z\in\D.
$$ 

We shall only need Lemma~\ref{Klemma} for the parameter value $p=\frac12$. 
In geometric language, the estimate in the lemma then says that 
the length of the radial segment $[0,1[$ in the metric 
(\ref{metric}) is less than or equal to $c_\alpha=c_{1/2,\alpha}$.
By a conformal mapping argument, %invariance 
we arrive at the proposition below. The notation ${\mathbf B}(z_0,r)$
stands for a metric disk of radius $r$ about the point $z_0\in\D$; also,
as before, the parameter $T=T(\omega)$ is the termination time for the 
Hele-Shaw flow emanating from the origin. 

\begin{prop}\label{HSdestimate}
Let ${\mathbf\Omega}=\langle\D,\omega\rangle$ be a 
simply connected Riemann surface with a metric 
$\diff{\mathbf s}^2=\omega\vert\diff z\vert^2$. 
Assume that the weight function $\omega$ is $C^2$-smooth and 
strictly positive in $\D$ and that the curvature condition 
$$
{\mathbf K}(z)+\alpha\,{\mathbf K}_\mathbb{H}(z)\leq0,\qquad z\in\D,
$$
is satisfied for some $\alpha$ with $0\leq\alpha<\pi/(8-2\pi)$. 
Let $c_\alpha=c_{1/2,\alpha}$ be as in Lemma~\ref{Klemma}.
%If the metric disk ${\mathbf B}(0,r)$ is precompact in ${\mathbf\Omega}$, 
%then, for $0<t<r^2/c_\alpha$, the Hele-Shaw domain $D(t)$ 
%is precompactly contained in ${\mathbf B}(0,r)$.
%In particular $D(t)\subset {\mathbf B}(0,c_\alpha\sqrt{t})$ provided that 
%${\mathbf B}(0,c_\alpha\sqrt{t})$ is precompact in $\Omega$.
If for some $t$, $0<t<T$, the Hele-Shaw domain $D(t)$ is simply connected and 
$\partial D(t)$ is a Jordan curve, then 
$D(t)\subset {\mathbf B}(0,c_\alpha\sqrt{t})$.
\end{prop}

\begin{proof}
%See the proof of~\cite[Proposition~8.2]{HS}.
Let $z_0\in D(t)$, and let $\varphi:\D\to D(t)$ be the conformal map 
such that $\varphi(0)=0$ and $\varphi^{-1}(z_0)\in[0,1[$.
Let $\nu$ be defined by (\ref{dfnnu}).
We shall apply Lemma~\ref{Klemma} to the weight $\nu$.
It is straightforward to see that $\nu$ satisfies the 
subharmonicity assumption in Lemma~\ref{Klemma}.
We proceed to check the reproducing property of $\nu$.

A classical result of Torsten Carleman~\cite[Section~1]{Carleman} 
asserts that the space of 
analytic polynomials is dense in the Bergman space $A^p(D)$ ($0<p<+\infty$) 
of all $p$-th power area integrable analytic functions in $D$ if 
the domain $D$ is simply connected with Jordan boundary. Let $h$ be a 
harmonic polynomial. 
In view of Carleman's theorem, we can find a sequence $\{h_j\}$ of functions 
harmonic in a neighbourhood of $\bar D(t)$ such that 
$h_j\to h\circ\varphi^{-1}$ in $L^1(D(t))$. 
By the mean value property of $D(t)$ (see~\cite[Theorem~2.3]{HS}), we have 
$$
t\,h_j(0)=\int_{D(t)}h_j(z)\,\omega(z)\,\diff\Sigma(z),
$$
and a passage to the limit yields that
$$
t\,h(0)=\int_{D(t)}h\circ\varphi^{-1}(z)\,\omega(z)\,\diff\Sigma(z).
$$
An obvious change of variables now gives the reproducing property of $\nu$.

An application of Lemma~\ref{Klemma} gives that 
$$\int_0^1\sqrt{\nu(r)}\,\diff r\leq c_\alpha.$$ 
This last inequality shows that the geodesic distance in the 
metric (\ref{metric}) from the origin to the point $z_0$ is 
less than $c_\alpha\sqrt t$.
We conclude that $D(t)\subset{\mathbf B}(0,c_\alpha\sqrt{t})$. 
\end{proof}

We remark that for $\alpha=0$ the above proposition 
recovers~\cite[Proposition~8.2]{HS}.
For a complete surface, we have the following corollary.

\begin{cor}
Let ${\mathbf\Omega}=\langle\D,\omega\rangle$ 
be as in Theorem~\ref{simplyconnected} and 
assume in addition that ${\mathbf\Omega}$ is complete. 
Then $T=+\infty$. This in turn implies that the Hele-Shaw flow covers 
all of ${\mathbf\Omega}$, that is, $\D=\bigcup_{t>0}D(t)$. 
\end{cor}

\begin{proof}
Since ${\mathbf\Omega}$ is complete, a famous theorem of Hopf, Rinow 
and de Rham asserts that the metric disks ${\mathbf B}(0,r)$, $0<r<+\infty$, 
are all precompact in ${\mathbf\Omega}=\langle\D,\omega\rangle$ 
(see~\cite[Theorem~I.10.3]{Helgason}). In view of Theorems 
\ref{simplyconnected} and \ref{realanalytic}, the Hele-Shaw domains $D(t)$ 
are simply connected with real-analytic Jordan curve boundaries for as long 
as $0<t<T$, that is, $D(t)\Subset\D$. Moreover, Proposition~\ref{HSdestimate} 
guarantees that the growth of the Hele-Shaw domains is essentially no faster 
than that of the metric disks. The conclusion is that $D(t)$ remains 
precompact in $\D$ for all $t$, $0<t<+\infty$, so that $T=+\infty$.
In addition, an obstacle problem argument shows that the Hele-Shaw flow 
covers all of ${\mathbf\Omega}$ (see~\cite[pp. 219-220]{HS}). 
\end{proof}

\section{The Hele-Shaw exponential mapping}\label{HSexpsection}

Just as in~\cite{HS}, we have the following existence result regarding 
the so called {\sl Hele-Shaw exponential mapping}.

\begin{thm}\label{HSexpmapping}
Let ${\mathbf\Omega}=\langle\D,\omega\rangle$ be as in 
Theorem~\ref{simplyconnected}.
Let $z_0\in{\mathbf\Omega}$.
Then there exists a real-analytic mapping 
$\Phi=\HSexp_{z_0}:\D(0,\sqrt{T})\to{\mathbf\Omega}$ such that
\begin{itemize}
\item $\Phi(0)=z_0$,
\item each ray $\{z\in\D(0,\sqrt{T})\setminus\{0\}:\ \arg(z)=\theta\}$ is 
mapped mapped 
onto a curve in ${\mathbf\Omega}$ which points in the same direction 
as the ray at $z_0$,
\item the map $\Phi$ maps each pair consisting of a concentric circle 
about the origin and a straight line passing through the origin  
onto a pair of orthogonal curves, and
\item for each $0<r<\sqrt{T}$, the domain $\Phi(\D(0,r))$ 
equals the Hele-Shaw domain ${\bf D}(z_0,r^2)$.   
\end{itemize}
The map $\Phi$ is uniquely determined in the class of $C^1$-mappings 
satisfying the above four properties. 
Furthermore, the map $\Phi$ has the asymptotics 
$$
\Phi(z)=z_0+\omega(0)^{-1/2}z+O(\vert z\vert^2)\qquad\text{as}\ 
\vert z\vert\to0.
$$
If the surface ${\mathbf\Omega}$ is complete, 
then $T=+\infty$, which in turn implies that the Hele-Shaw flow covers 
all of ${\mathbf\Omega}$, that is, ${\mathbf\Omega}$ is the union of all
the domains ${\bf D}(z_0,t)$ over all positive $t$. 
\end{thm}

\begin{proof}
Using properties of the weighted biharmonic Green function $\Gamma_1$ 
we have shown that the Hele-Shaw flow domain $D(t)$, $0<t<T$, 
is simply connected and has real-analytic Jordan boundary 
(see Section~\ref{wbhGfcn} and Section~\ref{simplyconnHSdomain}). 
Also, in Section~\ref{Korenblum} we showed that $T=+\infty$ if 
${\mathbf\Omega}$ is complete. Knowing this, the proof of 
Theorem~\ref{HSexpmapping} proceeds exactly as in~\cite{HS}. We omit 
the details. 
\end{proof}

\section{An example of a weakly hyperbolic surface}\label{exwhs}

In this section, we give an example of a complete real-analytic 
metric (\ref{metric}) satisfying (\ref{curvaturealpha}) 
for a given positive $\alpha$, which, for some $r$, $0<r<1$, has the circle 
$r\T$ as a geodesic curve. By an obvious pull-back with a conformal self-map 
of the unit disk, we obtain a metric of the kind referred to in the 
introduction.

Let $I$ be a bounded interval on the real line $\R$. The differential 
equation for a geodesic curve $\gamma:I\to\D$ 
in the metric (\ref{metric}) can be written
\begin{equation}\label{geodesiceq}
\gamma''(t)+\frac{1}{\omega(\gamma(t))}\,
\frac{\partial\omega}{\partial z}(\gamma(t))\,[\gamma'(t)]^2=0,
\qquad t\in I,
\end{equation}
where the curve $\gamma$ is to be parameterized proportionally to arc length, 
that is, the quantity 
$$\|\gamma'(t)\|_{\gamma(t)}^2=\omega(\gamma(t))\,\vert\gamma'(t)\vert^2$$ 
is to be constant. 

We now specialize to a radial weight:
$$\omega(z)=\omega_0\big(\vert z\vert^2\big),\qquad z\in\D.$$
Let $\gamma_r$ be the circle $r\T$ parameterized by 
$\gamma_r(t)=re^{\imag t}$ for real $t$, $-\pi<t\le\pi$. A straightforward 
computation shows that $\gamma_r$ satisfies (\ref{geodesiceq}) if and only if 
\begin{equation}\label{omega0eq}
\omega_0(r^2)+r^2\omega_0'(r^2)=0.
\end{equation}

Fix $\alpha$, $0<\alpha<+\infty$. We introduce a positive real parameter
$c$, and consider weights of the form 
$$
\omega(z)=\frac{c}{(1-\vert z\vert^2)^2}+(1-\vert z\vert^2)^{2\alpha},
\qquad z\in\D.
$$
By comparison with the Poincar\'e metric, it is clear that the metric 
(\ref{metric}) so obtained is complete. Since the function 
$$z\mapsto\frac{\omega(z)}{(1-\vert z\vert^2)^{2\alpha}}$$ 
is the sum of two logarithmically subharmonic functions, it is itself 
logarithmically subharmonic (see~\cite[Corollary~1.6.8]{HorSCV}). 
This shows that (\ref{curvaturealpha}) holds.

We now show that the equation (\ref{omega0eq}) has a solution $r$, $0<r<1$.
A computation shows that
$$
\omega_0(r)+r\omega_0'(r)=c\,\frac{3-r}{(1-r)^3}+
(1-r)^{2\alpha-1}\big(1-(1+2\alpha)r\big).
$$
Note that the second term on the right hand side is 
negative for $r>1/(1+2\alpha)$.
Thus, for $c$ small and positive, the function 
$$r\mapsto\omega_0(r)+r\omega_0'(r)$$
attains negative values in $]0,1[$.
Taking the endpoints $r=0,1$ also into account, we see that 
equation (\ref{omega0eq}) has at least two solutions $r$ in the interval 
$]0,1[$ for small positive values of $c$.

%%% References %%%

%\bigskip

%\noindent \textsc{H\aa{}kan Hedenmalm, Anders Olofsson,
%Department of Mathematics,  \\
%Royal Insti\-tute of Techno\-logy, SE--100 44 
%Stockholm, Sweden}

%\noindent\textsl{E-mail}: \texttt{haakanh@math.kth.se} (Hedenmalm), 
%\texttt{ao@math.kth.se} (Olofsson).

\end{document}